\documentclass[10pt]{article}
\usepackage{amsmath,amsfonts,latexsym,amssymb}
\usepackage[latin1]{inputenc}
\newtheorem{theorem}{Theorem}[subsection]
\newtheorem{lemma}[theorem]{Lemma}
\newtheorem{proposition}[theorem]{Proposition}
\newtheorem{coro}[theorem]{Corollary}
\newtheorem{conj}[theorem]{Conjecture}

\newcommand{\pr}{{\mathbb R\mathbb P}^1}

\newcommand{\inn}{_{n\in\mathbb N}}

\newcommand{\Rep}{\rm{Rep}}
\newcommand{\proof}{{\sc Proof : }}

\newcommand{\qed}{{\sc Q.e.d.}}
\newcommand{\maire}{{\rm MinArea}}
\newcommand{\grf}{\pi_1 (S)}
\newcommand{\bgrf}{\partial_\infty \pi_1 (S)}

\newcommand{\mapping}[4]
{
\left\{
\begin{array}{rcl}
#1 &\rightarrow& #2\\
#3 &\mapsto& #4 
\end{array}
\right.
}

\begin{document}
\newcommand{\auteur}{
\vskip 2truecm
\centerline{François Labourie}
\centerline{Topologie et Dynamique}
\centerline{Université Paris-Sud}
\centerline{F-91405 Orsay (Cedex)}}
\title{Cross Ratios, Anosov Representations and the Energy  Functional on Teichmüller Space.}

\author{François LABOURIE}
\maketitle
\section{Introduction}

Let $S$ be a closed surface of genus greater than one.  Let $\grf$ be its fundamental group.
We study in this article two types of representations of $\grf$ in Lie groups
\begin{itemize}
\item {\sc Hitchin representations of $\grf$ in $PSL(n,\mathbb R)$} 
Following \cite{FLA}, we define a {\em Fuchsian representation} of $\grf$ in $PSL(n,\mathbb R)$ to be a representation which factors through the irreducible representation of $PSL(2,\mathbb R)$ in  $PSL(n,\mathbb R)$ and a cocompact representation of $\grf$ in $PSL(2,\mathbb R)$. A {\em Hitchin representation} is a representation which may be deformed in a Fuchsian representation. In \cite{H}, Hitchin studies the moduli space of reductive ({\em i.e.} whose Zariski closure is a reductive group) Hitchin representations.  For another point of view on a related subject, see the work of Fock and Goncharov in \cite{FG}.
\item{\sc  Maximal symplectic representations of $\grf$} We recall  constructions and results from \cite{BIW2}.

The symmetric space $M$ associated to $PSp(2n,\mathbb R)$ is Hermitian and carries an invariant symplectic form $\omega$ as we now explain. Let $\partial_\theta$ be the generator of the centre of $K$ of $PSp(2n,\mathbb R)$ such that
$$
\exp(2\pi\partial_\theta)=1.
$$
The complex structure on $TM$ is given by
$$
A\mapsto [\partial_\theta,A].
$$
We normalize the Killing form $\langle\ ,\rangle$ so that $\Vert \partial_\theta\Vert =1 $.
Then the symplectic structure is given by
$$
\omega(X,Y)=\langle[X,Y],\partial_\theta\rangle.
$$
We now observe that the  homomorphism 
$$
\det : K\rightarrow S^1,
$$
 is a degree $n$ map when  restricted to the centre of $K$. It follows that $n.\omega$ is the curvature of a line bundle $L$.

Let $\rho$ be a representation of $\grf$ in  $PSp(2n,\mathbb R)$. Let  $f$ be a $\rho$-equivariant map from the universal cover of $S$ to $M$. Then  the $2$-form $f^*\omega$ is invariant under the action of  $\grf$.
The following number
$$\tau(\rho)=\frac{n}{2\pi}\int_S f^*\omega$$
is an integer independent of the choice of $f$, since it is the  first Chern class of the induced bundle from $L$ by $f$. This number $\tau(\rho)$ is called the {\em Toledo invariant} of $\rho$.  It is constant under continuous deformations of the representation. More importantly, the following inequality holds
$$
\vert \tau(\rho)\vert\leq n \vert\chi(S)\vert.
$$
For  $n=1$, this is the Milnor-Wood Inequality \cite{M}. This specific inequality for the symplectic group is due to V. Turaev \cite{T}. It has been extended to other Hermitian Symmetric spaces (see in particular
\cite{BIW1} and \cite{BIW2}). We shall restrict ourselves to $PSp(2n,\mathbb R)$ although the discussion in this paragraph extends to the general case as well.

Of specific interests are the {\it maximal representations} for  which the above inequality is an  equality. For $n=1$, it is due to W. Goldman that these representations are monodromies of hyperbolic structures \cite{WG1}, \cite{WG}.

These maximal representations have been extensively studied by Bradlow, Garc\'\i a-Prada, Gothen, Mundet i Riera (as well as Xia in a specific example)  (\cite{BGG}, \cite{BGG2}, \cite{BGG3}, \cite{GM}, \cite{Go}, \cite{Xia}) using Higgs bundle techniques on one hand and 
Burger,  Iozzi and Wienhard (\cite{BISU}, \cite{BIW1}, \cite{BIW2}, \cite{W}) using bounded cohomology techniques on the other hand. 

\end{itemize}
Both type of representations can be thought of generalisations of representations of $\grf$ in $PSL(2,\mathbb R)$ which are monodromies of hyperbolic structures. It is already known  that these two types of representations share some common properties. 
\begin{itemize}
\item They are Anosov as defined in \cite{FLA}. For Hitchin representations this is proved in \cite{FLA}, for maximal representations this is shown in  \cite{BILW} by Burger, Iozzi, Wienhard and the author.
\item They are reductive. For Hitchin representation see Proposition \ref{hitchred}, for maximal representations this is proved by Burger, Iozzi and Wienhard in \cite{BIW2}.
\item They are discrete, see \cite{FLA} for Hitchin representations, and the proof by Burger, Iozzi and Wienhard \cite{BIW2} for maximal representations. 

\end{itemize}

We show in this paper that they share some more common features. This will follow from the fact that they are {\em well displacing}. More precisely, let $\Gamma$ be a finitely generated subgroup of the isometry group of a metric space $X$, we say $\Gamma$ is {\em well displacing}, if given a set $S$ of generators of $\Gamma$, there exist positive constants $A$ and $B$ such that
$$
\forall \gamma\in\Gamma, \forall x\in X, d(x,\gamma(x))\geq A \inf_{\eta\in\Gamma}\Vert\eta\gamma\eta^{-1}\Vert_S -B.
$$
where $\Vert\gamma\Vert_S$ is the word length of $\gamma$ with respect to $S$. It is easy to check that cocompact groups are well displacing, as well as convex cocompact whenever $X$ is {\em Hadamard} (i.e. complete, non positively curved and simply connected).  If $\rho$ is a representation of $\Gamma$ in a Lie group $G$, we say it is {\em well displacing} if the  group $\rho(\grf)$ is  {\em well displacing} as groups of isometries of the associated symmetric space. Our first result is the following
\begin{theorem}\label{intwd}
Hitchin and maximal symplectic representations are well displacing. 
\end{theorem}
This follows by relating Hitchin representations and maximal
	symplectic representations to  cross ratios (cf Theorems  \ref{crhitch}
and \ref{crossmax}), a construction already  interesting in its own right. This theorem  also turns out to have some interesting consequences which we now list.

 First, we make a remark about quasi-isometries: it has  already been observed in \cite{BILW} the related fact that the orbit maps are quasiisometries for maximal symplectic representation. Moreover, Olivier Guichard has explained to me that, for a well displacing representation of a surface group (and a more general class of groups) in any isometry group,  all orbit maps are quasiisometries. 
 
 Now, we list in the next paragraph some applications of the  previous result.
 
\subsubsection{The Mapping Class Group acts properly}

Let ${\rm Hom}_H(\grf, SL(n,\mathbb R))$ be the space of Hitchin homomorphisms and 
$$
{\rm Rep}_H(\grf,SL(n,\mathbb R))={\rm Hom}_H(\grf, SL(n,\mathbb R))/SL(n,\mathbb R).
$$
Let ${\rm Hom}_T(\grf, PSp(2n,\mathbb R))$ be the space of maximal symplectic  homomorphisms  and 
$$
{\rm Rep}_T(\grf,PSp(2n,\mathbb R))={\rm Hom}_T(\grf, PSp(2n,\mathbb R))/PSp(2n,\mathbb R).
$$
Both spaces are union of connected components of space of representations and admit actions of the Mapping Class Group $\mathcal M
  (S)$. The following result is an almost obvious consequence of Theorem \ref{intwd}. 
\begin{theorem}\label{mcgint}
The Mapping Class Group $\mathcal M (S)$ acts properly on  the spaces ${\rm Rep}_H(\grf,SL(n,\mathbb R))$  and 
${\rm Rep}_T(\grf,PSp(2n,\mathbb R))$.
\end{theorem}

This was known in the case of $SL(3,R)$ combining Goldman result on the properness of the action of the Mapping Class Group on the moduli space of $\mathbb{RP}^2$ structures  \cite{WGo} and the identification of this latter space with Hitchin component by Choi and Goldman \cite{Ch-Go}.

\subsubsection{The energy is proper}

Let $\rho$ be a representation of $\grf$ in a semi-simple Lie group $G$ whose associated symmetric space is $M$. Let $M_\rho$ be the associated flat  $M$ bundle over $S$. Let $\Gamma (S,M_\rho)$ be the space of sections of $M_\rho$. If $J$ is a complex structure on $S$ and $f$ an element of $\Gamma(S,M_\rho)$ we define
$$
{\rm Energy}_J (f)=\int_S \langle df\wedge df\!\circ\! J \rangle.
$$
Associated to a  such a representation $\rho$, we define a function $e_\rho$, the  {\em energy}, on the Teichmüller space $\mathcal T (S)$ in the following way
$$
e_\rho (J) =\inf ({\rm Energy}_J (f), f \in \Gamma (S,M_\rho)).
$$
We give more precise definitions in Paragraph \ref{defenergy}. We prove

\begin{theorem}\label{propint}
If $\rho$ is a Hitchin representation or a maximal symplectic representation, then $e_\rho$ is a proper function on Teichmüller space.
\end{theorem}

It is classical that critical points of the energy relate to minimal surfaces. Indeed, as a corollary we obtain
\begin{coro}
Let $\rho$ be a Hitchin or maximal symplectic representation. There there exists a minimal branched immersion of $S$ in $M/\rho(\grf)$ which represents $\rho$ at the level of homotopy groups.
\end{coro}

We investigate this in our two cases and give some extra details.

\subsubsection{Minimal area and Toledo invariant}
We restrict in this paragraph to representations in $PSp(2n,\mathbb R)$. For such a representation $\rho$, we denote by 
$\tau(\rho)$ the Toledo invariant. We define the {\em minimal area of $\rho$} by 
$$
{\rm MinArea}(\rho)=\inf (e_\rho(J), J\in\mathcal T (S)).
$$
We  say  $\rho$ is {\em diagonal} if it factors through a cocompact representation in $PSL(2,\mathbb R)$ and the diagonal representation of $PSL(2,\mathbb R)$ in 
$$
\prod_{i=1}^{i=n} PSL(2,\mathbb R)\subset PSp(2n,\mathbb R).
$$ 
We prove
\begin{theorem}\label{toledminint}
For every $\rho$ we have $$
\frac{n}{2\pi}{\rm MinArea}(\rho)\geq \vert\tau(\rho)\vert.
$$
If furthermore $\rho$ is maximal and $$\frac{n}{2\pi}{\rm MinArea}(\rho)= \tau(\rho),$$
then $\rho$ is diagonal.
\end{theorem}

Oscar Garc{\'\i}a-Prada explained to me that this result should have a natural  Higgs bundle interpretation.

\subsubsection{The Hitchin map is surjective}

In his article  \cite{H}, N. Hitchin gives  explicit parametrisations of Hitchin
components. Namely, given a choice of a complex structure
$J$ over  a given compact surface $S$, he produces a homeomorphism
$$
H_J: \mathcal Q(2,J)\oplus\ldots\oplus\mathcal Q(n,J) \rightarrow{\rm Rep}_H(\grf,SL(n,\mathbb R)),
$$
where $\mathcal Q(p,J)$ denotes the space of holomorphic
$p$-differentials on the Riemann surface $(S,J)$.  The main
idea in the proof is  first to identify representations with
harmonic mappings as in K. Corlette's seminal paper
\cite{KC}, (see also \cite{D}, \cite{FL1}), second to use the fact  a harmonic
mapping $f$ taking values in a symmetric space gives rise to
holomorphic differentials in manner
similar to that in which  a connection gives rise to differential forms
in Chern-Weil theory (cf. Paragraph \ref{comCW}).

However one drawback of this construction is that $H_J$ depends on the choice of the complex structure $J$. In particular, it breaks the invariance by the Mapping Class Group and therefore this construction does not give information on the topological nature of  ${\rm Rep}_H(\grf,SL(n,\mathbb R))/\mathcal M
  (S)$. We explain now a more equivariant (with respect to the action of the Mapping Class Group) construction. Let 
$\mathcal E^{(n)}$ be the vector bundle over Teichmüller space whose fibre above the complex structure $J$ is
$$
\mathcal E^{(n)}_J=\mathcal Q(3,J)\oplus\ldots\oplus\mathcal Q(n,J).
$$
We oberve that the dimension of the total space of $\mathcal E^{(n)}$ is the same as that of ${\rm Rep}_H(\grf,SL(n,\mathbb R))$ since the dimension of the "missing" quadratic differentials in $
\mathcal E^{(n)}_J$ accounts for the dimension of Teichmüller space.
account for .
We now define the {\em Hitchin map} 
$$
H \mapping{\mathcal E^{(n)}}{{\rm Rep}_H(\grf,SL(n,\mathbb R)}{(J,\omega)}{H_J(\omega).}
$$
We are aware that this terminology is awkward since this Hitchin map  is some kind of an inverse of what is usually called the Hitchin fibration.
From Hitchin construction, it now  follows this map is equivariant with respect to the Mapping Class Group action. We prove

\begin{theorem}\label{hitchsurjint}
The Hitchin map is surjective.
\end{theorem}

Our strategy is to identify $\mathcal E^{(n)}$ with the moduli space of  equivariant minimal surfaces  in the associated symmetric space and to prove that there exists an equivariant minimal surface for every representation by tracking a critical point of the energy.

Our conjecture in \cite{FLA} is that the Hitchin map is a homeomorphism, which is a consequence of the following:
\begin{conj}
If $\rho$ is a Hitchin representation, then the  minimum of $e_\rho$ is non degenerate. 
\end{conj}

This conjecture is well known to be true for $n=2$. For
$n=3$, one can prove it using  ideas linking real
projective structures, affine spheres, Blaschke metrics as in
J. Loftin paper \cite{JL} or in the preprint \cite{FL2}.

If the last conjecture is true, then following our previous
discussion, we would obtain the following result, which
helps to understand the action of the Mapping Class Group
$\mathcal M (S)$ on Hitchin components
\begin{conj}
  The quotient $\Rep_{H}(\grf,PSL(n,\mathbb R))/\mathcal M
  (S)$ is homeomorphic to the total space of the vector
  bundle $E$ over the Riemann moduli space, whose fibre at a
  point $J$ is
  $$
  E_{J}=\mathcal Q(3,J)\oplus\ldots\oplus\mathcal Q(n,J).
  $$
\end{conj}   
Again, by the previous discussion this result is true for
$n=2$ and $n=3$.
\subsubsection{Comments and extensions}
The theory of Hitchin representations extends to all real split groups, and the theory of maximal representations extends to all isometry groups of Hermitian symmetric spaces. It t is quite natural to conjecture that the constructions of this article extend to these more general cases. The fact that these representations are (at least conjecturally) Anosov representations is certainly meaningful from this point of view. However, one cannot expect all the results here extend for all Anosov representations since one can construct Anosov representations which are not reductive. It  also remains a question to understand under which algebraic conditions Anosov representations are associated to cross ratios, which is a crucial argument in our paper. 

It would be interesting to study  the space $\mathcal S$ of equivariant minimal surfaces in the case of the other components than Hitchin component for maximal symplectic representations.   For the same reason as for Hitchin component, this map is surjective, but I have no idea of what is its structure.

I wish to thank L. Lemaire for help on minimal surfaces and harmonic mappings, F. Paulin for a very helpful comment on irreducible representations,  A. Iozzi for clarifications on Toledo invariant,  M. Burger, O. Garc{\'\i}a-Prada, O. Guichard  and A. Wienhard for numerous conversations and useful remarks on inconsistencies of the first draft, as well as  W. Goldman for his interest and his help on the writing up.

\subsubsection{Outline of the paper}
\begin{itemize}
\item{\sc  2. Cross ratio.}  We recall the basic definitions (cross ratios, periods)  and explain how cross ratios are related to flows. This helps to control the growth of the periods (Proposition \ref{spec})
\item{\sc  3. Representations and cross ratios.} We study specifically Hitchin and maximal symplectic representations. We explain that they generate curves in generalised Grassmannians,  are reductive and relate to cross ratios.
\item{\sc 4.  Energy, minimal area.} We recall basic results about existence of equivariant harmonic mappings in symmetric spaces, as well as classical definitions and results concerning minimal surfaces, energy and Teichmüller space.
\item{\sc 5.  Well displacing representations.} We introduce the notion of well displacing representations, show that Hitchin and maximal symplectic representations are well displacing using our previous result on cross ratios. Then we prove the main properties of well displacing representations and obtain Theorems \ref{mcgint} and \ref{propint}.
As a consequence, we deduce the existence of equivariant minimal surfaces for our two main examples.
\item{\sc 6. Toledo invariant and  minimal area.} We investigate the relation between Toledo invariant and minimal area and prove Theorem \ref{toledminint} .
\item{\sc 7. The Hitchin map.} We prove Theorem \ref{hitchsurjint}.
\end{itemize}
\section{Cross Ratio}
Let  $\bgrf$ be the boundary at infinity of $\grf$. We recall that it is a circle equipped with an action of $\grf$. This action can be characterised by the following two properties
\begin{itemize}
\item every orbit is dense,
\item every non trivial element of $\grf$ has exactly two fixed points: one attractive, one repulsive.
\end{itemize} 
The boundary at infinity $\bgrf$  is identified with the boundary of the Poincaré disk model when one chooses a uniformisation of the universal cover of the surface. Let
$$\bgrf^{4*}=\{(x,y,z,t)\in\bgrf^{4}, x\not=t , \hbox{ and } y\not= z\}.$$
A {\em strict cross ratio} on $\bgrf$ is a $\grf$-invariant   Hölder function $b$ on $\bgrf^{4*}$ with values in $\mathbb R$ which satisfies the following rules
\label{birrules}
\begin{eqnarray}
b(x,y,z,t)&=&b(z,t,x,y) \label{bir100}\\
b(x,y,z,t)&=&0\ \ \Leftrightarrow x=y \hbox{ or } z=t\\
b(x,y,z,t)&=&{b(x,y,z,w)}{b(x,w,z,t)}\label{bir11}\\
b(x,y,z,t)&=&b(x,y,w,t)b(w,y,z,t)\label{bir11bis}\\
b(x,y,z,t)&=&1 \ \Leftrightarrow x=z \hbox{ or } y=t\label{bir12}.
\end{eqnarray}
The classical cross ratio on $\pr$, given in projective coordinates by 
\begin{equation*}
	b(x,y,z,t) = \frac{(y-x)(t-z)}{(y-z)(t-x)},
\end{equation*} is an example of a strict cross ratio.

For more information and examples, compare with \cite{CR}, \cite{led}.
\vskip 1 truecm
\subsection{Periods}
Let $b$ be a cross ratio $b$ and $\gamma$ be a nontrivial element in $\grf$. The {\em period} $l_b(\gamma)$ is defined as follows.  Let  $\gamma^{+}$ (resp. $\gamma^{-}$) the attracting (resp.  repelling) fixed point of $\gamma$ on $\partial_\infty\pi_1 (S)$.  Let $y$ be an element of $\partial_\infty \pi_1 (S)$. We define
\begin{eqnarray}
l_b(\gamma,y)=\log\vert b(\gamma^-,\gamma y,\gamma^+,y)\vert.\label{period}
\end{eqnarray}
It is immediate to check that $l_b(\gamma)=l_b(\gamma,y)$ does not depend on $y$. Moreover, by Equation( \ref{bir100}), $l_b(\gamma)=l_b(\gamma^{-1})$.

The next proposition compares periods with length of geodesics.
\begin{proposition}\label{spec}
For every non trivial $\gamma$ in $\grf$, let $\lambda(\gamma)$ be the length of the closed geodesic associated to $\gamma$ for a given hyperbolic metric.
Let $b$ be a strict cross ratio. Then there exists a  positive constant $A$, depending only on the cross ratio and the choice of the hyperbolic metric,  such that
$$
\forall\gamma\in\grf,\ \frac{1}{A}\lambda(\gamma)\leq l_b(\gamma)\leq A \lambda(\gamma).
$$
\end{proposition}
In the next paragraph, we define {\em compatible flows} on $\bgrf^{3*}$ and study their periodic orbits. Then, we show that every cross ratio is associated to a compatible flow and finally conclude.
\subsubsection{Compatible flow}
A flow $\phi_t$ on $\bgrf^{3+}$ is {\em compatible} if
\begin{itemize}
\item it is $\grf$ invariant,
\item it has no fixed points,
\item $\phi_t(x,z,y)=(x,u,y)$ with $(x,z,u,y)$ cyclically ordered. 
\end{itemize}
For the  convenience of notations, associated to a compatible flow $\phi_t$ we define the map 
$$ \hat\phi_t: \bgrf^{3*}\rightarrow \bgrf$$ such that
$(x,\hat\phi_t(x,z,y),y)=\phi_t(x,z,y)$.
If $\gamma$ is a non trivial element of $\grf$ we define the {\em period} of $\gamma$ to be the number $l_\phi(\gamma)$ such that
$$\hat\phi_{l_\phi(\gamma)}(\gamma^-,y,\gamma^+)=\gamma(y)$$.

We first prove
\begin{proposition}\label{speccompa}
Let $\phi_t$ and $\psi_t$ be two compatible flows. Then there exists a constant $K$ such that
$$
\forall \gamma\in\grf\setminus\{id\}, l_\psi(\gamma)\leq K l_\phi (\gamma).
$$
\end{proposition}

\proof There exists a continuous $\grf$-invariant positive function $T$ defined on $\bgrf^{3*}$ such that
$$
\hat\psi_1(x,y,z)=\hat\phi_{T(x,y,z)}(x,y,z).
$$
Since $\bgrf^{3*}/\grf$ is compact, it follows there exists a constant  $A$ such that, 
$$
\forall (x,y,z)\in\bgrf^{3*},\  0< T(x,y,z)\leq A.
$$
Then we  have
$$
l_\psi(\gamma)\leq A l_\phi(\gamma) +A.
$$
By compactness of $\bgrf^{3*}/\grf$, there exists a constant $l$ such that
$$
\forall \gamma\in\grf\setminus\{id\},\  l_\phi (\gamma)\geq l >0.
$$
Therefore it follows that
$$
l_\psi(\gamma)\leq (A +A/l) l_\phi(\gamma).
$$
\qed

\subsubsection{Proof of Proposition \ref{spec}}
\begin{proposition}\label{homeo}
Let $b$ be a strict cross ratio. Let $x$, $y$ be two distinct elements of $\bgrf$. Let $I$ be one of the connected component of $\bgrf\setminus\{x,y\}$. Let $z$ be an element of $I$.
Then the map $\phi$
$$
\mapping{I}{\mathbb R }{t}{\phi(t)=\log (b(x,t,y,z)),}
$$ is a homeomorphism.
\end{proposition}

\proof  We first prove that $\phi$ is injective. Suppose that $\phi(s)=\phi(t)$. This implies that
$$
b(x,t,y,s)=\frac{b(x,t,y,z)}{b(x,s,y,z)}=e^{\phi(t)-\phi(s)}=1.
$$
Hence $s=t$.

It follows that $\phi(I)=]\alpha,\beta[$. We now prove that $\beta=+\infty$. Assume on the contrary that $\beta<\infty$. It follows that
$$
\lim_{t\rightarrow y}{\log (b(x,t,y,z))}=\beta <\infty
$$
Let us  choose an auxiliary compatible flow $\psi_t$ on $\bgrf^{3*}$. It follows that
\begin{eqnarray}
\lim_{t\rightarrow y}{b(x,t,y,\hat\psi_1 (x,t,y))}=1\label{eq:homeo}
\end{eqnarray}
Let now choose a sequence $\{t_n\} \inn$ converging to $y$. Since $\bgrf^{3^*}/\grf$ is compact, there exists a sequence $\{\gamma_n\}\inn$ of elements of $\grf$ such that
$$
\lim_{n\rightarrow\infty}{(\gamma_n (x),\gamma_n(y),\gamma_n(t_n))}=(X,Y,T)\in \bgrf^{3*}.
$$
If follows from Assertion (\ref{eq:homeo}), that
\begin{eqnarray*}
b(X,T,Y,\hat\psi_1(X,T,Y))&=&\lim_{n\rightarrow\infty}{b(\gamma_n(x),\gamma_n(t_n),\gamma_n(y),\hat\psi_1 (\gamma_n(x),\gamma_n(t_n),\gamma_n(y))}\\
&=&\lim_{n\rightarrow\infty}{b(\gamma_n(x),\gamma_n(t_n),\gamma_n(y),\gamma_n(\hat\psi_1 (x,t_n,y))}\\
&=&\lim_{n\rightarrow\infty}{b(x,t_n,y,\hat\psi_1 (x,t_n,y))}\\
&=&1.
\end{eqnarray*}
Hence $T=\hat\psi_1(X,T,Y)$ and the contradiction since $\psi_1$ has no fixed points. A similar argument yields $\alpha=-\infty$. Therefore $\phi$ is a homeomorphism. \qed

\begin{proposition}\label{exiscompa} There exists a compatible  flow $\phi_t$ on $\bgrf^{3+}$ such that  $\log(b(x,z,y,\hat\phi_t(x,z,y)))=t$.
\end{proposition}
\proof By Proposition \ref{homeo}, $\hat\phi_t$ is well defined. It is a flow by the relation
$$
b(x,y,z,t).b(x,t,z,u)=b(x,y,z,u).
$$
\qed.

Now, Proposition \ref{spec} is a consequence of Propositions \ref{exiscompa} and \ref{speccompa}

\section{Representations and Cross Ratio}

We explain that our two favorite classes of representations of $\grf$ are associated to crossratios whose periods can be computed from the holonomy. We also recall that these representations are reductive.

\subsection{Hitchin representations}

\label{defhitch}
Following \cite{FLA}, we define a {\em Fuchsian representation} of $\grf$ in $PSL(n,\mathbb R)$ to be a representation which factors via the irreducible representation of $PSL(2,\mathbb R)$ in $PSL(n,\mathbb R)$ and a cocompact representation of $\grf$ in $PSL(2,\mathbb R)$. A {\em Hitchin representation} is a representation that can be deformed into a Fuchsian representation $\pi_1(S)\to PSL(n,R)$. We furthermore say it is {\em reductive} if the Zariski closure of its image is a reductive group. Actually, we shall later on show that every Hitchin representation is reductive.  In \cite{H}, Hitchin studies the moduli space of reductive Hitchin representations.

In \cite{FLA}, we show

\begin{theorem}\label{split}
Let $\rho$ be a  reductive Hitchin representation. Let $\gamma$ be a non trivial element of $\grf$. Then $\rho(\gamma)$ is $\mathbb R$-split.
\end{theorem}

\subsubsection{Reductivity}
Recall  that we say a subgroup (or a representation) of $SL(n,\mathbb R)$ is  {\em irreducible} if it does not preserve any proper subspace in the standard representation of $SL(n,\mathbb R)$.

We show

\begin{proposition}\label{hitchred}
Every Hitchin representation is irreducible.
\end{proposition}

\proof In \cite{FLA}, Lemma 10.1, we explain how it follows from elementary observations on Higgs bundle  that every {\em reductive} Hitchin representation is irreducible. We now explain how to get rid of the "reductive" hypothesis. 

To prove that every Hitchin representation is irreducible, it suffices to show that the set of reductive Hitchin representations is closed. Let $\rho$ be a  limit of Hitchin representations. Let $G$ be the Zariski closure of $\rho(\grf)$. Let $N$ be the nilradical of $G$. Let $R=G/N$ be the reductive part of $G$. We identify $R$ with a subgroup of $G$ so that $G=N\rtimes R$. Let $\pi$ be the projection from $G$ to $R$. 
We first observe that $\pi\circ\rho$ is also a limit of reductive representations. 

Indeed, there exists an element $h$ in the centraliser of $R$ such 
$$
\forall u\in N, \lim_{n\rightarrow\infty}h^{-n}uh^{n}=1.
$$
It follows that 
$$
\pi\circ\rho=\lim_{n\rightarrow\infty}h^{-n}\rho h^{n}.
$$

It follows that $\pi\circ\rho$ is also a Hitchin representation hence irreducible.  Assume $N$ is non trivial. Since $N$ is unipotent, the set of vector fixed by $N$ is a proper subspace of the standard representation. This set is fixed by $R$, therefore $\pi\circ\rho$ is not irreducible. Hence we obtain the contradiction

We owe this argument to F. Paulin. \qed

\subsubsection{Cross Ratio}

In \cite{CR}, we associate to every Hitchin representation a cross ratio. More precisely, we show

\begin{theorem}\label{crhitch}
For every Hitchin representation $\rho$, there exists a cross ratio $b$ such that for every, a nontrivial element  $\gamma$ of $\grf$ then the period of $\gamma$ is given by
\begin{eqnarray}
l_b(\gamma)=\log(\vert\frac{\lambda_{max}(\rho(\gamma))}{\lambda_{min}(\rho(\gamma))}\vert),\label{eq:crhitch}
\end{eqnarray}
where $\lambda_{max}(\rho(\gamma))$ and $\lambda_{min}(\rho(\gamma))$ are  the  eigenvalues  of respectively maximum and minimum modulus of the element $\rho(\gamma)$.
\end{theorem}

We also prove in \cite{CR} a reciprocal of this statement, which amongst its ingredients uses O. Guichard's work \cite{OG}.

\subsection{Symplectic Anosov Structures}

In \cite{BILW}, we studied maximal representations of surface groups in $PSp(2n,\mathbb R)$. Let's start with some definitions and notations. 
\subsubsection{Toledo invariant}\label{deftoled}
We recall  constructions and results from \cite{BIW2}. We explained in the introduction the normalisation we choosed for the metric and the construction of the Toledo invariant.
The following inequality holds
$$
\vert \tau(\rho)\vert\leq n\vert \chi(S)\vert.
$$
For $n=1$, this is the Milnor-Wood Inequality \cite{M}. This specific inequality for the symplectic group is due to V. Turaev \cite{T}. It has been extended to other Hermitian symmetric spaces (see in particular \cite{D-To}, \cite{To1}, \cite{To2},
\cite{BIW1} and \cite{BIW2}). We shall restrict ourselves to $PSp(2n,\mathbb R)$ although the discussion extends to the general case as well.

Of specific interests are the {\it maximal representations} for which the above inequality is an  equality.

\subsubsection{Reductivity}
Using bounded cohomology techniques, Burger,  Iozzi and Wienhard prove in \cite{BIW2}
\begin{theorem}\label{maxred}
Every maximal representation is reductive.
\end{theorem}

\subsubsection{Cross Ratio and Maximal Representations}
For every $A$ in $PSp(2n,\mathbb R)$. Let $\{\lambda_i\}_{1\leq i\leq 2n}$ be the eigenvalues (with multiplicities) of $A$ ordered so that
$$
\vert \lambda_1\vert\leq\vert \lambda_2\vert\ldots \leq \vert \lambda_{2n}\vert.
$$
We define
$$ 
c(A)=\prod_{i=n+1}^{i=2n}\vert\lambda_i\vert.
$$
The main result of this paragraph is the following.
\begin{theorem}\label{crossmax}
Let $\rho$ be a maximal symplectic representation. Then there exists a cross ratio $b$ such that 
\begin{eqnarray*}
l_b(\rho(\gamma))=2\log c(\rho(\gamma)).
\end{eqnarray*}
\end{theorem}
We prove this Theorem in Paragraph \ref{proofcrossmax}

\subsubsection{Positivity}
Let $\mathcal L(E)$ be the Grassmanian of Lagrangian spaces in a vector space $E$ equipped with a symplectic form $\omega$. 
We say a triple of Lagrangian spaces $(F,G,L)$ is {\em positive} if
\begin{itemize}
\item $F\oplus L=E$
\item for every non zero vector $u$ in $G$ $\omega(u_f,u_l)>0$, where $u=u_f+u_l$ with $u_f\in F$ and $u_l\in L$.
\end{itemize}
We finally say an oriented curve $\xi : S^1\rightarrow L(E)$ is {\em positive} if for every oriented triple $(x,y,z)$ in $S^1$, $(\xi(x),\xi(y),\xi(z))$ is positive.  We prove in \cite{BILW} with Burger, Iozzi and Wienhard,  the following result
\begin{theorem}\label{bilw}
Let $\rho$ be a maximal symplectic representation of $\pi_1(S)$. Then there exists a positive equivariant curve $\xi:\bgrf\rightarrow\mathcal L(E)$.

Furthermore, $\xi(\gamma^+)$ (respectively $\xi(\gamma^-)$) is generated by the eigenvectors of $\rho(\gamma)$ corresponding to the eigenvalues of absolute value greater than 1 (resp. smaller than 1).
\end{theorem}

\subsubsection{Cross Ratio}\label{crosslag}
Let $(L_1,L_2,L_3,L_4)$ be a quadruple of Lagrangian spaces in a symplectic space of dimension $2n$. We suppose that $L_4$ is transverse to $L_1$ and $L_2$ transverse to $L_3$. Let $(l^1,l^2,l^3,l^4)$ be basis of $(L,E,F,G)$ respectively. For every pair $(a,b)\in \{1,2,3,4\}^2$, we consider the $n\times n$  matrix
$$
A_{l^a,l^b}=(\omega(l^a_i,l^b_j)).
$$
We observe  that for every endomorphism $g$ of  $L^a$ whose matrix in the basis $l^a$ is $G$,
\begin{equation}
A_{g(l^a),l^b}=G. A_{l^a,l^b}.\label{eq:crosslag}
\end{equation}
Similarly
$$
A_{l^a,l^b}= -A^t_{l^b,l^a}.
$$
We now define
$$
B(l_1,l_2,l_3,l_4)=\frac{\det(A_{l^1,l^2}).\det(A_{l^3,l^4})}{\det(A_{l^1,l^4}).\det(A_{l^3,l^2})}.
$$
By Assertion (\ref{eq:crosslag}), 
$B(l_1,l_2,l_3,l_4)$ only depends on
$(L_1,L_2,L_3,L_4)$.
Hence, we define
$$
B(L_1,L_2,L_3,L_4)=B(l_1,l_2,l_3,l_4).
$$
Moreover, 
\begin{eqnarray*}
B(L_1,L_2,L_3,L_4)B(L_1,L_4,L_3,L_5)&=&B(L_1,L_2,L_3,L_5)\\
B(L_1,L_2,L_3,L_4)&=&B(L_2,L_1,L_4,L_3).
\end{eqnarray*}
Finally, if $(L,U,V)$ are generic,
\begin{eqnarray*}
B(L,U,L,V)=1\\
B(L,L,U,V)=0.
\end{eqnarray*}
\subsubsection{Positivity and Cross Ratio}
\begin{proposition}\label{crosspos}
If $(E,F_1,G)$ is positive as well as $(E,F_2,G)$ then
$$
B(E,F_1,G,F_2)>0
$$
Furthermore, if  $(F_1,F_2,G)$ is positive, then 
$$
B(E,F_1,G,F_2)>1
$$
\end{proposition}
\proof We assume that $(E,F_1,G)$ is positive as well as $(E,F_2,G)$. Let $p$ be the projection on  $G$ along $E$. Let $q_i$  be the quadratic form on $F_i$.
defined by 
$$
q_i(u)=\omega(p(u),u).
$$
Since $(E,F_i,G)$ is positive, $q_i$ is  positive definite. By simultaneous orthogonalisation, we can choose an orthogonal  basis $f^i$  for  $F_i$ for $q_i$ such that $p(f^1)=p(f^2)=g$. 
Let finally $e^i=(1-p)(f^i)=f^i-g$ be  the corresponding bases of $E$. Therefore, we have
$$
A_{e^i,f^j}= A_{e^i,g}=-A_{g,f^i}.
$$
We also have 
$$
\omega(e^j_i,g_k)=q (f^j_i,f^j_k)=\omega(g_i,f^j_k).
$$
Therefore if 
$$
\lambda_i=\frac{q(f^2_i,f^2_i)}{q(f^1_i,f^1_i)},
$$
we have
$$
e^2_i=\lambda_i e^1_i.
$$
It follows that 

\begin{eqnarray*}
B(E,F_1,G,F_2)&=&\frac{\det(A_{e^1,f^1}).\det(A_{g,f^2})}{\det(A_{e^1,f^2}).\det(A_{g,f^1})}\cr
&=&\frac{\det(A_{e^1,g}).\det(-A_{e^2,g})}{\det(A_{e^1,g}).\det(-A_{e^1,g})}\cr
&=&\prod_i \lambda_i  >0.
\end{eqnarray*}
Finally, assume that $(F_1,F_2,G)$ is positive.
Let $\pi$ be the projection on $G$ along $F_1$. Recall that
$$
f^2_i=e^2_i + g_i={\lambda_i}e^1_i +g_i= {\lambda_i}f^1_i +(1-{\lambda_i})g_i.
$$
It follows that
\begin{eqnarray*}
\omega\big(\pi(f^2_i),(1-\pi)(f^2_i)\big)&=&\omega\big((1-{\lambda_i})g_i ,\lambda_i f^1_i\big)\\
&=&(1-\lambda_i) \lambda_i\omega(g_i,f^1_i)\\
&=&\lambda_i (\lambda_i-1) q_1(f^1_i,f^1_i).
\end{eqnarray*}
Hence the positivity of $(F_1,F_2,G)$ implies that $\lambda_i > 1$ and therefore
$$
B(E,F_1,G,F_2)=\prod_i \lambda_i  >1.
$$
\qed

Finally, we have.

\begin{proposition}\label{periodcrosslag}
Let $S$ be a symplectic automorphism preserving two transverse Lagrangian spaces $E$ and $F$. Then for every Lagrangian space $G$ we have
$$
B(E,G,F,S(G))=\frac{\det(S\vert_E)}{\det(S\vert_F)}=\det(S\vert_E)^2.
$$
\end{proposition}
\proof Let $S$ be a  symplectic transformation. Let $e$ be a basis of a space $K$ invariant by $S$. We have
$$
\det (A_{e,S(l)})=\det(A_{S^{-1}e,l})=\frac{\det(A_{e,l})}{\det(S\vert_K)}.
$$
Hence the formula follows.
\qed
\subsubsection{Cross ratio and maximal representations}\label{proofcrossmax}
We now prove Theorem \ref{crossmax}. Let $\rho$ be a maximal symplectic representation. By Theorem \ref{bilw}, let $\xi$ be the associated positive map in $\mathcal L(\mathbb R^{2n})$. By Proposition \ref{crosspos}, the following formula which uses the notations of Paragraph \ref{crosslag} defines a strict crossratio on $\bgrf$
$$
b(x,y,z,t)=B(\xi(x),\xi(y),\xi(z),\xi(t)).
$$
Furthermore, by Proposition \ref{periodcrosslag}, we have
$$
l_b(\gamma)=2\log\det (\rho(\gamma)\vert_{\xi(\gamma^+)}.
$$
The statement follows.

\section{Energy, minimal area}\label{defenergy}

Let $M$ be a Hadamard ({\em i.e}complete,  non positively curved and  simply connected) manifold with metric $g_M$.  Let $\rho$ be a representation of $\grf$ in the group ${\rm Iso} (M)$ of isometries of $M$.  Let $M_\rho$ be the associated $M$ bundle over $S$. We introduce some definitions which are more or less standard.

We define $\mathcal F_\rho$ to be the space of $\rho$-equivariant mapping of $S$ in $M$:
$$
\mathcal F_\rho=\{f:\widetilde S\rightarrow M/ f\circ\gamma=\rho(\gamma)\circ f\}.
$$
Another way to think of  $\mathcal F_\rho$ is to identify it with  $\Gamma(S,M_\rho)$ the space of sections of $M_\rho$.
 
\subsection{Energy}
Let $f$ be an element of $\mathcal F_\rho$. Let $J$ be  a complex structure on $S$. We define the following exterior differential 2-form on $S$
$$
\langle df\wedge df\!\circ\! J\rangle (x,y)=\langle T\! f(y),T\! f(Jx)\rangle -\langle T\! f(x),T\! f(Jy)\rangle.$$
We notice that $\langle df\wedge df\!\circ\! J\rangle $ is a $\grf$-invariant. We then define the {\em energy} of $f$ with respect to $J$ to be 
$$
E(J,f)=\int_S \langle df\wedge df\!\circ\! J\rangle.
$$

Finally, associated to $\rho$, we define the {\em energy} functional on Teichmüller space  ${\mathcal T} (S)$ associated to $\rho$ by 
$$
e_\rho: \mapping{{\mathcal T} (S)}{\mathbb R}{J}{e_\rho(J)=\inf_{f\in\mathcal F_\rho}(E(J,f))}.
$$ 

\subsection{Harmonic mappings}
By definition, a harmonic mapping is a critical point of the energy. Whenever $\rho$ is reductive, the existence of a $\rho$-equivariant mapping is guaranteed
by Corlette's Theorem in the context of symmetric spaces \cite{KC}. Note that  \cite{FL1} gives an alternative simpler proof which works in the general context of  Hadamard manifolds:

\begin{theorem}{\sc [Corlette]}. If $\rho$ is reductive  there exists a $\rho$-equivariant harmonic mapping $f$ from $S$ to $M=G/K$ the symmetric space associated to $G$. Furthermore this mapping is unique  up to an isometry of $G/K$ and minimizes the energy. 
\end{theorem}

Combining Corlette's Theorem with Propositions \ref{maxred}, \ref{hitchred}, we deduce the following result

\begin{proposition}\label{exisharm}
Let $\rho$ be a Hitchin or maximal symplectic representation. Let $J$ be a complex structure on  $S$. Then there exists a unique (up to isometries) $\rho$-equivariant harmonic mapping $f_{\rho,J}$. Moreover
 $$
 E(J,f_{\rho,J})=\inf_{f\in\mathcal F_\rho} E(f,J).
 $$
 \end{proposition}
\subsection{Area}
Let $f$ be an element of $\mathcal F_\rho$. Let $R(f)$ be the open set of points $x$ in $S$ for which $T_xf$ is  injective. We observe that the induced bilinear form $f^*(g_M)$ is invariant by $\grf$ and actually defines a metric on $R(f)$. We then define the {\em area} of $f$ to be the area of $R(f)$ for this metric.
$$
{\rm Area}(f)={\rm area}_{f^*(g_M)}(R(f)).
$$
It is well known that $f$ satisfies
$$
{\rm Area}(f)\leq E(J,f),
$$
with equality only if $f$ is conformal with respect to $J$.

We finally define the {\em minimal area of $\rho$} to be 
$$
\maire (\rho)=\inf_{J\in\mathcal T(S)}{e_\rho(J)}.
$$
Obviously, one has 
$$
\maire (\rho)=\inf_{ f \in\mathcal (F)}{\rm Area}(f).
$$
We finally recall  classical results by Sacks-Uhlenbeck \cite{SU1}\cite{SU2} and Schoen-Yau \cite{SY}.

\begin{theorem}
\label {SU}{\sc{[Sacks-Uhlenbeck],  [Schoen-Yau] }} Let $J$ be a point in Teichmüller space which is a critical point of  the  energy. Let $f_J$ be a mapping such that 
$$
E(J,f)={\rm MinArea}(\rho).
$$
Then $f$ is harmonic and conformal.
\end{theorem}

\section{Well displacing representations}
\subsection{Displacement function}\label{defwd}
Let $\gamma$ be an isometry of a metric space $M$. We recall that the {\em displacement} of $\gamma$ is
$$
d(\gamma)=\inf_{x\in M}(x,\gamma(x)).
$$
If $M$ is the Cayley graph of a group $\Gamma$ with set of generators $S$ and word length $\Vert\ \Vert_S$, then
$$
d(\gamma)=\inf_{\eta}\Vert\eta\gamma\eta^{-1}\Vert_S.
$$
The displacement function is explicit in the case of $M=SL(n,\mathbb R)/SO(n,\mathbb R)$.
Let $A$ be an element in $SL(n,\mathbb R)$. Let $\{\lambda_i\}_{1\leq i\leq n}$ be the eigenvalues of $A$, then
\begin{eqnarray}
d(A)= \log(\frac{1}{n}\sum_{i=1}^{i=n} \vert \lambda^2_i\vert).\label{eq:expli}
\end{eqnarray}
We  say a representation of $\Gamma$ in $Iso(M)$ is {\em well displacing} if for every set $S$ of generators of $\Gamma$ there exist positive constants $A$ and $B$ such that
$$
\forall \gamma\in\grf, \ d(\rho(\gamma)) \geq A .\inf_{\eta}\Vert\eta\gamma\eta^{-1}\Vert_S -B.
$$
where $\Vert \gamma\Vert_S$ is the word length of the element $\gamma$ of $\grf$. In other words, the displacement in the space $M$ is roughly controlled by the displacement in the Cayley graph.
Alternatively in the case $\Gamma=\grf$ and if  $\lambda(\gamma)$ is the length of the closed geodesic associated to the non trivial element of $\gamma$  , $\rho$ is well displacing if for every hyperbolic metric on $S$, there exist positive constants $A$ and $B$ such that 
$$
\forall \gamma\in\grf, \ d(\rho(\gamma)) \geq A .\lambda(\gamma) -B.
$$
where  $\lambda(\gamma)$ is the length of the closed geodesic associated to the non trivial element of $\gamma$.
\subsection{Examples of well displacing representations}
It is immediate to check that cocompact groups are well displacing, as well as convex cocompact whenever $X$ is Hadamard.
As examples of well displacing representation, and using our results on cross ratio, we have
\begin{proposition}\label{hitchwd}
Every Hitchin representation is well displacing. Every maximal symplectic representation is well displacing.
\end{proposition}

\proof Let $\rho$ be a Hitchin representation. Let $b$ be the associated cross ratio by Theorem \ref{crhitch}. Let $\{\lambda_i\}_{1\leq i\leq n}$ be the eigenvalues of the element $\rho(\gamma)$. We suppose $\lambda_1=\lambda_{min}$ and $\lambda_n=\lambda_{max}$. Combining Equations \ref{eq:crhitch} and \ref{eq:expli}, we obtain :
\begin{eqnarray*}
d(\rho(\gamma))+\log(n)&=&\log(\sum_{i=1}^{i=n} \lambda^2_i)\\
&\geq& \frac{2}{n} \log \vert \lambda_n\vert^n \\
&\geq&  \frac{2}{n}\log \big(\vert \lambda_n\vert \prod_{i\not=1}\vert\lambda_i\vert\big)=\frac{2}{n}\log \frac{\vert \lambda_n\vert}{\vert\lambda_1\vert}=\frac{2}{n}l_b(\gamma).
\end{eqnarray*}
By Theorem \ref{spec}, there exists a positive constants $A$ and $B$, such that
$$
l_b(\gamma)\geq A. \lambda(\gamma) -B ,
$$
 where $\lambda(\gamma)$ is the length of the geodesic associated to $\gamma$ in some auxiliary hyperbolic metric.
It follows that every Hitchin representation is well displacing.

For maximal symplectic representations, we have a very similar argument.  Let $\rho$ be such a representation. Let $\hat b$ be the cross ratio associated by Proposition \ref{crossmax}.
Since the injection $i$ of $PSp(2n,\mathbb R)$ in $PSl(2n,\mathbb R)$ gives rise to a totally geodesic embedding of the corresponding symmetric spaces, it suffices to show $i\circ \rho$ is well displacing.  Let $\{\lambda_i\}_{1\leq i\leq 2n}$ be the eigenvalues, with multiplicities, of the element $\rho(\gamma)$. We  order them  so that
$$
\vert \lambda_1\vert\leq\vert \lambda_2\vert\ldots \leq \vert \lambda_{2n}\vert.
$$
Then
\begin{eqnarray*}
d(\rho(\gamma))+\log(2n)&=&\log(\sum_{i=1}^{i=2n} \lambda^2_i)\\
&\geq& \frac{1}{n} \log \vert \lambda_{2n}\vert^{2n} \\
&\geq&  \frac{1}{n}\log \big(\vert \prod_{i=n+1}^{i=2n}\vert\lambda_i\vert^2\big)=\frac{1}{n}l_b(\gamma).
\end{eqnarray*}
The result follows. \qed

The main result of this section is the following.
\begin{theorem}\label{wdenergyproper}
Let $\rho$ be a well displacing representation in the isometry group of a Hadamard manifold. Then 
the energy functional $e_\rho$, defined  from $\mathcal T (S)$ to $\mathbb R$, is proper.
\end{theorem}

In particular, we recover as a corollary that the energy is proper for convex cocompact representations as was shown by W. Goldman and R. Wentworth in \cite{Go-W}.

We come back now to our favorite examples. We deduce the following result.

\begin{coro}\label{exismin}
Let $\rho$ be a Hitchin or maximal representation. Then there exists a complex structure $J_0$ on $S$, a conformal harmonic $\rho$-equivariant mapping $f$ such that
$$
{\rm Area}(f)=E(J_0,f)={\rm MinArea}(\rho).
$$
\end{coro}

\proof Indeed By Proposition \ref{hitchwd} and the previous Theorem, there exists a  complex structure $J_0$ on $S$ which achieves the minimum of the energy. By Proposition \ref{exisharm}, there exists a $\rho$-harmonic mapping $f$, such that
$$
E(J_0,f)=e_\rho(J_0)={\rm MinArea}(\rho).
$$
We conclude by Theorem \ref{SU} of Sachs-Uhlenbeck, Schoen-Yau.\qed

\subsubsection{The intersection is proper}

Let $g$ and $g_0$ be two hyperbolic metrics on $S$. Let $US$ and $U_0 S$ be the associated unit tangent bundle with geodesic flows $\phi_t$ and $\phi^0_t$, generated by $X$ and $X_0$. Let $\mu$ and $\mu_0$ be the  corresponding Liouville measure. We know that the geodesic flows are orbit  conjugate. In other words, there exist a map $F$ from $US$ to $US_0$ and a  positive function $\psi_{(g,g_0)}$ on  $US$ such that $F$ is differentiable along $X$ and $DF(\psi_{g,g_0}X)=X_0$. We define the intersection of $g$ and $g_0$ to be
$$
{\rm inter}(g ,g_0)=\int_{US} \psi_{(g,g_0)} d\mu.
$$
The following Proposition is classical. For the sake of completeness since we could not find a good reference for it, we include a sketchy proof. For a less down to earth point of view and extra information  on intersection,  we refer to Curt Mc Mullen notes \cite{CMC} and Francis Bonahon original article \cite{FB}.
 
\begin{proposition}\label{interproper}
The function $g\mapsto {\rm inter}(g ,g_0)$ is proper on $\mathcal T(S)$.
\end{proposition}
\proof 
By definition, the intersection of two closed curved $\gamma_{1}$ and $\gamma_{2}$ is 
$$
{\rm inter(\gamma_{1},\gamma_{2})}=\inf(\sharp(c_{1}\cap c_{2}), \hbox{ such that } c_{i} \hbox { is homotopic to } \gamma_{i}).
$$
If $\gamma_{i}$ are geodesics for a negatively curved metric $g$, then 
$$
{\rm inter(\gamma_{1},\gamma_{2})}=\sharp(\gamma_{1}\cap \gamma_{2}).
$$
Let $\lambda_{g}(\gamma)$ be the length of the geodesic $\gamma$ with respect to the metric $g$.
Using a tubular neighbourhood along $\eta$, we obtain that there exists a constant $C(\eta,g)$ just depending on $\eta$ and the metric $g$, such that for every closed curve $\gamma$, if 
\begin{equation}
 {\rm inter}(\eta,\gamma)\leq C(\eta,g)\lambda_{g}(\gamma).\label{intereq2}
\end{equation}
Let now $g$ and $g_{0}$ be two hyperbolic metrics. Let 
 $$
 \mathcal G_L=\{\gamma / \lambda_{g}(\gamma)\leq L\}.
 $$
 We recall that for every function $f$, the Liouville measure satisfies
 $$
 \int_{US} fd\mu=\lim_{L\rightarrow\infty}\bigg(\frac{1}{\sharp( \mathcal G_L)}\sum_{\gamma\in\mathcal G_L}\frac{\int_\gamma fdt}{\lambda_{g}(\gamma)}\bigg).
 $$
 Hence
 $$
{\rm inter}(g ,g_0)=\lim_{L\rightarrow\infty}\bigg(\frac{1}{\sharp( \mathcal G_L)}\sum_{\gamma\in\mathcal G_L}\frac{\lambda_{g_{0}}(\gamma)}{\lambda_{g}(\gamma)}\bigg).
$$
 Furthermore, let $\eta$ be a closed curve, we also have 
\begin{equation}
 \lambda_g(\eta)=\lim_{L\rightarrow\infty}\bigg(\frac{1}{\sharp( \mathcal G_L)}\sum_{\gamma\in\mathcal G_L}\frac{\rm inter(\eta,\gamma)}{\lambda_{g}(\gamma)}\bigg).\label{intereq1}
\end{equation}
Combining Equation (\ref{intereq1}) and Inequality (\ref{intereq2}), we obtain that for every closed curve $\eta$,
 $$
 \lambda_{g}(\eta)\leq C(\eta,g_{0}) {\rm inter}(g ,g_0).
 $$
 Since we can find a finite set of closed curves $J$ such that the function
 $$
 \lambda_{J}: g\rightarrow \sum_{\eta\in J}\lambda_{g}(\eta) .
 $$
 is proper on Teichmüller space, the statement follows by the following inequality
 $$
 \lambda_{J}(g) \leq \big(\sum_{\eta\in J}C(\eta,g_{0})\big) {\rm inter}(g ,g_0).
 $$
 \qed

\subsubsection{The energy is proper}

We adapt here  a beautiful argument by C. Croke and A. Fathi \cite{CF}. Let $\rho$ be a well displacing representation. We use the notation of the previous paragraph. Let $J$ be a complex structure on $S$. Let $g$ be the associated hyperbolic metric.  Let $g_{0}$ be a fixed hyperbolic metric. We now prove, adapting \cite{CF}, that there exist a constant $K$ independent on $J$ such that
$$
e_\rho(J)\geq K ({\rm inter}(g ,g_0))^2.
$$
Let $f$ an element of $\mathcal F\vert_\rho$.  We consider the function
$$
h:\mapping{US}{\mathbb R}{u}{\Vert T\! f(u) \Vert}.
$$
We have
$$
E(J,f)=\int_{US}h^2d\mu.
$$
By Cauchy-Schwarz, 
$$
E(J,f)\geq (\int_{US}h\ d\mu)^2.
$$
Let $\gamma$ be a closed orbit of the geodesic flow. Let $\lambda_g(\gamma)$ be the length of $\gamma$ with respect to $g$. We denote also by $\gamma$ the corresponding element in the fundamental group. Let $c$ be the curve
$$
\mapping{[0,\lambda_g(\gamma)]}{M}{t}{f(\gamma(t)).}$$
Then 
$$
\int_\gamma h\  dt={\rm length}(c)\geq d(c(0),\rho(\gamma)(c(0))\geq A \lambda_{g_{0}}(\gamma) -B=A \int_\gamma \psi_{(g,g_0)} dt -B .
$$
We recall  that for every function  $F$ :
$$
\int_{US}F\  d\mu=\lim_{L\rightarrow\infty}\bigg(\frac{1}{\sharp( \mathcal G_L)}\sum_{\gamma\in\mathcal G_L}\frac{\int_\gamma F\ dt}{\lambda_{g}(\gamma)}\bigg).
 $$
It follows that for every $f$
\begin{eqnarray*}
\sqrt{E(J,f)}&\geq& \int_{US}h\ \ d\mu=\lim_{L\rightarrow\infty}\bigg(\frac{1}{\sharp( \mathcal G_L)}\sum_{\gamma\in\mathcal G_L}\frac{\int_\gamma h dt}{\lambda_{g}(\gamma)}\bigg).
\\
&\geq&A\lim_{L\rightarrow\infty}\bigg(\frac{1}{\sharp( \mathcal G_L)}\sum_{\gamma\in\mathcal G_L}\frac{\int_\gamma \psi_{(g,g_0)}dt}{\lambda_{g}(\gamma)}\bigg)
-B\lim_{n\rightarrow\infty}\bigg(\frac{1}{\sharp( \mathcal G_L)}\sum_{\gamma\in\mathcal G_L}\frac{1}{\lambda_{g}(\gamma)}\bigg).\\&\geq&
A \int_{US} \psi_{(g,g_0)} d\mu=A \ {\rm inter}(g ,g_0).
\end{eqnarray*}
Hence
$$
e_\rho(J)\geq A^2 ({\rm inter}(g ,g_0))^2.
$$
Finally, we know  by Proposition \ref{interproper} that the function $g\rightarrow {\rm inter}(g ,g_0)$ is proper. Hence $e_\rho$ itself is proper.

\subsection{Mapping Class Group and well displacing representations}
Let $X$ be a topological space. Let $x$ and $y$ be two elements of $X$, we write $x\mathcal R y$ if every neighbourhood of $x$ meets any neighbourhood of $y$. We consider  $\sim$ the equivalence relation generated by $\mathcal R$. We observe that $X^{red}=X/\sim$ is Hausdorff and, moreover, every continuous map from $X$ to a Hausdorff space factors through $X^{red}$.

Let  $Hom^*(\grf, {\rm Iso}(M))$ be the space of well displacing homomorphisms. Let
$$
{\rm Rep}^{red}(\grf,{\rm Iso}(M))=(Hom^*(\grf, {\rm Iso}(M)/{\rm iso(M)})^{red}.$$
 We have

\begin{proposition}
The Mapping Class Group ${\mathcal M}(S)$ acts properly on the space  ${\rm Rep}^{red}(\grf,{\rm Iso}(M))$.
\end{proposition}

\proof  Let ${\mathbb R^{\grf}}$ be the space of maps from $\grf$ to $\mathbb R$. We equip it with the product topology. As before, for every $\gamma$ in $\grf$, denote by $\lambda(\gamma)$  the length of the closed geodesic associated to $\gamma$ in a fixed hyperbolic metric. Let 
$$
{\mathbb R}^{\grf}_{hyp}=\{l\in{\mathbb R^{\grf}}/ \exists A, B >0,  \ \ \forall \gamma \in\grf\setminus\{id\}, \ \  l(\gamma)\geq  A \lambda(\gamma) -B.\}
$$

We observe that for every divergent sequence $g_n$ of elements of  ${\mathcal M}(S)$, there exists an element $\gamma$ such that
$$
\lim_{n\rightarrow\infty} \lambda(g_n(\gamma))=\infty.
$$
It follows that ${\mathcal M}(S)$ acts properly on ${\mathbb R}^{\grf}_{hyp}$.
Let 
$$
S:\mapping{{\rm Hom}^*(\grf, {\rm Iso}(M))}{\mathbb R^{\grf}}{\rho}{\{d(\rho(\gamma))\}_{\gamma\in\grf}.}
$$
It follows from the definition, that $S$ takes values in  ${\mathbb R}^{\grf}_{hyp}$. The Proposition follows. \qed

\section{Toledo invariant, minimal area and immersed minimal surfaces}

\begin{theorem} Let $\rho$ be a representation in  $PSp(2n,\mathbb R)$. Then
$$
\frac{n}{2\pi}{\rm MinArea}(\rho)\geq \tau(\rho).
$$
Moreover, if  $\tau(\rho)$ is maximal and ${\rm MinArea}(\rho)= \tau(\rho)$, then $\rho$ is diagonal.
\end{theorem}
We begin with a lemma.
\begin{lemma}\label{omk}
Let $M$ be an Hermitian symmetric space. Let $\omega$ be its symplectic form. Let  $(X,Y)$ be an  orthonormal pair of tangent vectors and $K$ be the sectional curvature of the plane generated by $(X,Y)$, then 
\begin{eqnarray}
\vert\omega(X,Y)\vert\leq \sqrt{- K}.\label{eq:omk}
\end{eqnarray}
Moreover the equality occurs exactly whenever $[X,Y]$ generates the centre of $K$. 
\end{lemma}
\proof Write $M=G/K$. Let $\partial_\theta$ be the generator of the centre of $K$ such that
$$
exp(2\pi\partial_\theta)=1.
$$
We normalize the Killing form $\langle,\rangle$ so that $\Vert \partial_\theta\Vert =1 $.
Then the symplectic structure is given by
$$
\omega(X,Y)=\langle[X,Y],\partial_\theta\rangle\leq \sqrt {-\langle[X,Y],[X,Y]\rangle}=\sqrt {-K}.
$$
\qed

We can now prove our Theorem.

\proof The first point (i) is immediate. Let $\omega$ be the Kähler form on the associated symmetric space $X$. Let $(u,v)$ be an orthonormal system in $T_x X$. Then, 
$\omega(u,v)\leq 1$ with equality if and only if the plane generated by $(u,v)$ is complex. It follows that for every $\rho$-equivariant mapping we have
$$
\frac{n}{2\pi}{\rm Area}(f)\geq\frac{n}{2\pi}\int_Sf^*(\omega)=\tau(\rho).
$$
If $f$ is an immersion, and the equality holds, then 
$f$ is an holomorphic map when $\tilde S$ is equipped with the induced complex structure.  

Assume now that $\rho$ is maximal. According to Theorem \ref{exismin}, there exist a complex structure $J$ on $S$ and a conformal harmonic mapping $f$ such that
$$
{\rm Area}(f)={\rm MinArea}(\rho).
$$

Let $f$ be as above. We know by results of \cite{Gu}, that $f$ is a branch minimal immersion. Let $x_1,\ldots,x_n$ be the branched points of order $k_1,\ldots,k_n$. Let $\hat S=S\setminus\{x_1,\ldots,x_n\}$.
Let's now denote by $\kappa$ the curvature of the metric $f^*g_X$ on $\hat S$.  
Notice that 
$$
\frac{1}{2\pi}\int_{\hat S}\kappa d\mu -\sum_i (k_i-1)=\chi(S)
$$
Let $\kappa_f$ be the sectional curvature of the 2-plane $T\! f(T\hat S)$. Let $B$ the second fundamental form of  $f$. By the Gauss equation
$$
\kappa=\kappa_f -\Vert B\Vert \leq \kappa_f.
$$
Finally, assume now that $\tau(\rho)=\frac{1}{2\pi}{\rm MinArea}(\rho)$. Let $\mu$ be the measure of area of the metric $f^*g_X$. We have by Inequality (\ref{eq:omk})
\begin{eqnarray*}
2\pi \vert \tau(\rho)\vert &=&n\vert \int_{\hat S} f^*\omega\vert\\
&\leq& n\int_{\hat S} \sqrt{- \kappa_f}d\mu\\
&\leq& n \sqrt{{\rm Area}(f)\int_{\hat S} - \kappa_f d\mu}\\
&\leq& \sqrt{2n\pi\tau(\rho)}\sqrt{\int_{\hat S} - \kappa_f d\mu}
\end{eqnarray*}
It follows that
\begin{eqnarray*}
\frac{1}{n}\vert\tau(\rho)\vert&\leq& -\frac{1}{2\pi}\int_{\hat S}\kappa_f d\mu\\\
&\leq& -\frac{1}{2\pi}\int_{\hat S}\kappa d\mu-\frac{1}{2\pi}\int_{\hat S}\Vert B\Vert d\mu\\
&\leq&-\chi(S)-\sum_i (k_i-1)-\frac{1}{2\pi}\int_{\hat S}\Vert B\Vert d\mu\\
&\leq &\frac{1}{n}\vert\tau(\rho)\vert-\sum_i n(k_i-1)-\frac{1}{2\pi}\int_{\hat S}\Vert B\Vert d\mu.
\end{eqnarray*}

It first follows that for all $i$, $k_i=1$. In orther words, $f$ is an immersion. Moreover $B=0$. This means that $f$  is totally geodesic. It follows that $f$ is associated to an embedding of $PSL(2,\mathbb R)$ in $PSp(2n)$ whose Lie algebra contains the centre of $K$ by the equality case in Lemma \ref{omk}.  Hence $\rho$ is diagonal.\qed

\section{The Hitchin map}

\subsection{Representations and holomorphic differentials}

We first recall how representations of $\grf$ in a semi-simple Lie group $G$ give rise to holomorphic differentials, by a construction quite similar to Chern-Weil theory. The construction that we now describe associate to every reductive representation of $\grf$ in $G$, to every complex structure $J$ on $S$ to every $ad$-invariant polynomial $q$ of degree $n$ on the Lie algebra of $G$ a holomorphic $n$-ic differential on $S$.

By Corlette's Theorem (\cite{KC} and \cite{FL1} for an alternative simpler proof), there exists a $\rho$-equivariant harmonic mapping $f$ from $S$ to $M$ the symmetric space associated to $G$. Moreover this mapping is unique  up to an isometry of $G/K$.  We define
 $${\rm Hom}_{red}(\grf,G)=\{\hbox{ reductive homomorphisms } \grf\rightarrow G\},$$
 and 
 $${\rm Rep}_{red}(\grf,G)={\rm Hom}_{red}(\grf,G)/G.$$

\subsubsection{Harmonic maps on surfaces}
We begin by  a standard observation on harmonic maps on surfaces. 

Let $f$ be a map to a space $M$. Let $T\! f$ be the tangent map of $f$. We consider $T\! f$ as a $1$-form on $S$ with values in the pullback bundle by $f$ of $TM$
$$
T\! f\in\Omega^1(S,f^*TM).
$$
Let now 
$$
\Omega^1_{\mathbb C}(S,f^*TM\otimes_{\mathbb R}\mathbb C),
$$
be  the space of complex $1$-form on $S$ with value in the complexified vector bundle $f^*TM_{\mathbb C}=f^*TM\otimes_{\mathbb R}\mathbb C$. We consider the  natural map
$$
\mapping{\Omega^1(S,f^*TM)}{\Omega^1_{\mathbb C}(S,f^*TM_{\mathbb C})}{\omega}{\omega_{\mathbb C},}
$$
defined by
$$
\omega_{\mathbb C}(u)=\omega(u)-i\omega (Ju).
$$
We say an element $\beta$ of $
\Omega_{\mathbb C}^1(S,f^*TM_{\mathbb C})
$ is {\em holomorphic} if
$$
\nabla_{Ju}\beta = i\nabla_u \beta.
$$
We have the following classical observation
\begin{proposition}
$$
f {\hbox{ is harmonic}} \Leftrightarrow T\! f_{\mathbb C} {\hbox{ is holomorphic}}.
$$
\end{proposition}
\proof Indeed, $f$ is harmonic if and only if for every $X$
$$
\nabla_{X}T\! f (X)+\nabla_{JX}T\! f (JX)=0.
$$
Since $\nabla_{X}T\! f (Y)$ is symmetric in $X$ and $Y$, the above condition is equivalent to
$$
\forall X,Y\in TS,\ \ \ \nabla_{X}T\! f (Y)+\nabla_{JX}T\! f (JY)=0.
$$
wich turns to be equivalent to 
$$
\forall X,Y\in TS, \ \ \ \nabla_{JX}T\! f(Y)-\nabla_{X}T\! f(JY)=0,
$$
On the other hand
$$
(\nabla_{JX}T\! f_{\mathbb C} -i\nabla_{X}T\! f_{\mathbb C}) (Y)=(\nabla_{JX}T\! f(Y)-\nabla_{X}T\! f(JY))-i(\nabla_{X}T\! f (Y)+\nabla_{JX}T\! f (JY)).
$$
The statement follows. \qed

\subsubsection{Commutative Chern-Weil Theory}\label{comCW}
Let $M$ be a Riemannian manifold. 
Let $p$ be a parallel section of $(T^*M)^{\otimes n}$. We denote by $p_{\mathbb C}$ the associated parallel  section of  $(TM_{\mathbb C}^*)^{\otimes n}$.  Let $f$ be a map from a Riemann surface $S$ to $M$. Then we have the following standard observation.
\begin{proposition}
Let $
\beta\in\Omega^1_{\mathbb C}(S,f^*TM_{\mathbb C})$. Suppose that $\beta$  is  holomorphic.
Then $p(\beta,\beta, \ldots, \beta)$ is a holomorphic differential of degree $n$.
\end{proposition}

As a specific example of this construction, we can take $p=g$, the Riemannian metric on $M$. We have the following immediate result

\begin{proposition}\label{gmin}
Let $f$ be a harmonic map map from a surface $S$ to $M$. Then $g_{\mathbb C}(T_{\mathbb C}f,T_{\mathbb C}f)=0$ if and only if $f$ is minimal.
\end{proposition}
\proof Indeed, $g_{\mathbb C}(T_{\mathbb C}f,T_{\mathbb C}f)=0$ if and only if $f$ is conformal, hence minimal. \qed

When $M$ is the symmetric space associated to $G={\rm Iso(M)}$,  then every $G$-invariant symmetric multilinear form $P$ on the Lie algebra $\mathcal G$ of $G$ 
gives naturally rise to a parrallel polynomial function $f_P$ on $M=G/K$.  Indeed such a $P$ give naturally rise to a $G$-invariant tensor $f_p$ on $G/K$. But any tensor on a symmetric space invariant under isometries is parallel.

Combining the above constructions and Corlette's theorem, it follows that for every complex structure $J$ on $S$, every symmetric $G$-invariant multilinear form $P$ of degree $n$ on $\mathcal G$, we have a map 
$$
F_{P,J}:\mapping{{\rm Rep}_{red}(\grf,G)}{{\mathcal Q}(n,J)}{\rho}{f_P(T\! f_{\mathbb C},\ldots,T\! f_{\mathbb C}),}
$$
where $f$ is a $\rho$-equivariant harmonic mapping from $S$ to $G/K$, given by Coreltte's Theorem. We denote here, as in the introduction, by ${\mathcal Q}(n,J)$ the space of holomorphic $n$-differentials on on $S$ equipped with the complex structure $J$. The uniqueness part of Corlette's result shows that $F_{P,J}$ is well defined.

When $G=SL(n,\mathbb R)$, let $p_n$ be the symmetric polynomial of degree $n$. Notice that $F_{p_{2}}$ is the metric on $G/K$. We define the  map 
$$
\xi_J=\bigoplus_{k=2}^{k=n}F_{p_k,J}.
$$
We can now state Hitchin Theorem \cite{H}
\begin{theorem}
{\sc{ [Hitchin]}} The  map $\xi_J$ is a homeomorphism from the space of Hitchin representations
${\rm Rep}_H(\grf,PSL(n,\mathbb R))$ to 
$$
\mathcal Q(2,J)\oplus\ldots\oplus\mathcal Q(n,J) .
$$
\end{theorem}
As in the introduction, we define the {\em Hitchin map} 
$$
H \mapping{\mathcal E^{(n)}}{{\rm Rep}_H(\grf,SL(n,\mathbb R)}{(J,\omega)}{\xi^{-1}_J(\omega).}
$$
This map is equivariant with respect to the Mapping Class Group action.

The following result is now immediate.

\begin{theorem}
The Hitchin map is surjective.
\end{theorem}

\proof Let $\rho$ be a Hitchin representation. By Corollary \ref{exismin} there exists a complex structure $J$ on $S$ and a $\rho$ equivariant conformal harmonic mapping $f$. It follows by \ref{gmin}, that $F_{p_2,J}(\rho)=0$. This shows the Hitchin map is surjective. \qed
\subsubsection{Normal bundle to the space of Fuchsian representations}
We conclude by a partial result. When $\rho$ is Fuchsian in $PSL(n,\mathbb R)$, the corresponding energy is the same as the energy for $\rho$ in $PSL(2,\mathbb R)$. Hence it has a unique strict minimum. Therefore the same is true for representations which are closed to be Fuchsian. It follows the Hitchin map is a diffeomorphism on a small neighbourhood of the zero section. This implies that the normal bundle of the space of Fuchsian representations in Hitchin component can be identified - equivariantly with respect to the action of the Mapping Class Group - with ${\mathcal E^{(n)}}$.

\end{document}